\documentstyle [oldlfont,12pt,amsfonts,amssymb,oldlfont]{article}

\newtheorem{lemm}{Lemma}      
   
\newtheorem{theo}{Theorem}   
\newtheorem{coro}{Corollary}
\newtheorem{defi}{Definition}

\newenvironment{proo}{\noindent {\bf Proof:}}{$\Box$ \vspace{7mm}}

\newmathalphabet*{\scri}{eus}{m}{n}
\newcommand{\labl}[1]{\label {#1}}

\def\rest{\mathord{\restriction}}

\newcommand{\text}{\mbox}

\newcommand {\dom} {\text{\rm{dom}}}
\newcommand {\rng} {\text{\rm{rng}}}
\newcommand {\fld} {\text{\rm{fld}}}

\newcommand {\IDM} {\text{\rm{IDM}}}

\def \rest {\mathord{\restriction}}

\begin {document}

\title{Identities on Cardinals Less Than $\aleph _{\omega}$}
\author { M. Gilchrist  \and S. Shelah}  
\maketitle

\section {Introduction}

Let \footnote{S. Shelah partially supported by
a research grant from the basic research fund of the Israel Academy of Science;
Pul. Nu. 491.}
 $ \kappa $ be an uncountable cardinal and the edges of a complete graph with $ \kappa $ vertices be colored with $ \aleph _{0} $ colors.  For $ \kappa > 2 ^{\aleph _{0}} $ the Erd\H{o}s-Rado theorem implies that there is an infinite monochromatic subgraph.  However, if $ \kappa \leq 2 ^{\aleph _{0}} $, then it may be impossible to find a monochromatic triangle.  This paper is concerned with the latter situation.  We consider the types of colorings of finite subgraphs that must occur when the edges of the complete graph on $ \kappa \leq 2^{\aleph _{0}} $ vertices are colored with $ \aleph _{0}$ colors.  In particular, we are concerned with the case $ \aleph _{1} \leq \kappa \leq \aleph _{\omega }$.  

The study of these color patterns (known as identities) has a history that involves the existence of compactness theorems for two cardinal models \cite{Schmerl}.  When the graph being colored has size $ \aleph _{1} $, the identities that must occur have been classified by Shelah \cite {Shappendix}.  If the graph has size greater than or equal to $ \aleph _{\omega} $ the identities have also been classified in \cite {Shtwocard}. The number of 
colors is fixed at $\aleph _{0}$ as it is the natural place to start and the 
results here can be generalized to situations where more colors are used.

There is one difference that we now make explicit.  When countably many
colors are used we can define the following coloring of the complete
graph on $2 ^{\aleph _{0}} $ vertices.  First consider the branches in the complete binary tree of height $\omega$ 
to be vertices of a complete graph. The edge $ \{ \eta, \nu \} \in [\,
^{\omega }2]^{2} $ is given the color $ \eta \cap \nu $.  (In the notation that follows this coloring induces the identities in $\IDM$ and no others.)
When  $  \kappa > \aleph _{0} $ colors are used we may not be able to define
an analogous coloring of the complete graph on $ 2 ^{ \kappa } $ vertices
since it is known that there exist models of set theory in which no tree 
has $ \aleph _{1} $ nodes and $ 2^{\aleph _{1}}$ branches \cite
{Baalmost}.

The first section of this paper deals with the definitions that are necessary to describe identities.  In the remainder of the paper we use forcing to construct a model $V$ of set theory in which the set of identities which must occur when $ \kappa = \aleph _{m+1}$, strictly includes the set of identities which must occur when $ \kappa = \aleph _{m}$.  This is done by showing that a certain identity is omitted when $ \kappa = \aleph _{m}$ but is induced when $ \kappa = \aleph _{m+1}$.
\section{Preliminaries}

An $\omega$-{\em coloring} is function $f:[B]^{2} \longrightarrow
\omega$ where $B$ is a set of ordinals ordered in the usual way. 
The set $B$ is the {\em field} of $f$ and is denoted fld$(f)$.

\begin{defi} Let $f,\,g$ be 
  $\omega$-colorings.  We say that $f$ {\em realizes} the coloring $g$ if
  there is a one-one map $ k : \fld (g) \longrightarrow \fld (f) $
  such that for all $ \{ x,y\} ,\{ u,v \} \in \dom (g)$
  $$ 
    f(\{ k(x), k(y) \} ) \not=  f(\{ k(u), k(v) \} ) \Rightarrow g(\{
    x,y \}) \neq  g(\{ u,v \} ).
  $$  
\end{defi}
We write $f \simeq g$ if $f$ realizes $g$ 
and $g$ realizes $f$.  It should be clear that $ \simeq$ induces an
equivalence relation on the class of $\omega$-colorings.  We call the
equivalence classes {\em identities}.  The collection of all identities is
denoted $ \text{ID}$.

\begin{defi} Let $f,\,g$ be 
$\omega$-colorings.  We say that $f$ {\em V-realizes} the coloring $g$ if
there is an order-preserving map $ k : \fld(g) \longrightarrow
\fld (f) $ such that for all $ \{ x,y\} ,\{ u,v \} \in \dom (g)$
$$ f(\{ k(x), k(y) \} ) \not=  f(\{ k(u), k(v) \} ) \Rightarrow g(\{ x,y \} )
 \not=  g(\{ u,v \} ).$$  We write $f \simeq _{V} g$ if $f$ V-realizes $g$
and $g$ V-realizes $f$.  Note that $ \simeq_{V}$ induces an
equivalence relation on the class of $\omega$-colorings.  We call the
equivalence classes {\em V-identities}.  The collection of all V-identities
is denoted $ \text{ID}_{V}$.
\end{defi}

For both types of realization we will call the map $k:\fld (g)
\longrightarrow \fld (f) $ an {\em embedding}.  In the above
definition the $V$ refers to vertices.  In the situations that
follow, $B$ will be a cardinal less than or equal to $ \aleph _{\omega
} $ ordered in the usual way as a set of ordinals.  In the following
we will speak of $\omega$-colorings realizing (rather than
V-realizing) other $\omega$-colorings whenever the context makes the
type of realization clear.  If $ f,\,g,\,h,\,l$ are $ \omega
$-colorings, with $ f \simeq g $ and $ h \simeq l$, then $f$ realizes
$h$ if and only if $g$ realizes $l$.  Thus without risk of confusion
we may speak of identities realizing colorings and of identities
realizing other identities. The same is true of V-identities.  If $I$
and $J$ are identities we call $J$ a {\em subidentity} of $I$ if $I$ realizes $J$.  The notion of
sub-V-identity is similarly defined.  We say that an identity $I$ is
of {\em size r } if $| \fld (f) | = r $ for some (all) $f \in I $.  In
the following we will consider only identities of finite size.

\medskip
 
An identity can be regarded as a finite structure $ \langle A,E
\rangle $, where $E$ is an equivalence relation on $[A]^{2}$, see
\cite {Shappendix}. The correspondence is given by: $A= \fld (f)$, and
$\{x,y \} \simeq_{E} \{ u,v \}$ if and only if $f(\{ x,y \} ) = f( \{
u,v \} )$. 

A V-identity can be regarded as a finite structure $\langle A,E,<_{A}
\rangle$, where $E$ is an equivalence relation on $[A]^{2}$ and
$<_{A}$ is a linear ordering of $A$.  The correspondence is given by: 
$A=\fld(f)$, $\{x,y \} \simeq_{E} \{ u,v \}$ if and only if $f(\{ x,y \} )
= f( \{ u,v \} )$, and $x< y$ if and only if $ x <_{A} y $.

We remark that the notion of subidentity does not correspond to that
of substructure.  To simplify the language in what follows we find it
convenient to abuse terminology by referring to $ \omega$-colorings
and structures $ \langle A,E \rangle $ as identities rather than as
representatives of identities.  Similarly for V-identities.  In each
case the intended meaning should be clear.

\begin{defi}  
  Let $ f:[B]^{2} \longrightarrow \omega $ be an $\omega$-coloring and
$I= \langle A,E, < _{A} \rangle $ be a structure corresponding to a
V-identity.  Let $k : A \longrightarrow B$ be an order-preserving
map such that for all $\{ x,y \} , \{ v,w \} \in [A]^{2}$, 
$$ 
  \{ k(x), k(y) \} \not\simeq_{E}  \{ k ( v), k(w) \} \Rightarrow
  f(\{ x,y\}) \not= f( \{ v,w \}) .  
$$
Then $k$ is called an {\em embedding} of $I$ into $ f$.  A similar
definition is given when $I$ is an identity.
\end{defi}

\begin{defi} Let $ D$ be a set of ordinals and $ f: [B ] ^{2} \longrightarrow \omega$ be an $\omega$-coloring.
\begin{enumerate}
\item ${\cal I }(f)\, \,({\cal I }_{V}(f) )\,\,$ is the collection of (V-) identities realized by $f$.
\item ${\cal I }(D) = \bigcap \{ {\cal I } (g) \,|\, g:  [D ] ^{2} \longrightarrow \omega \}, \, \,{\cal I }_{V}(D) = \bigcap \{ {\cal I }_{V} (g) \,|\, g:  [D ] ^{2} \longrightarrow
\omega
 \}$.
\end{enumerate}
\end{defi}
\medskip
If $J=\langle B, F \rangle$ is an identity and $ A \subset B$ we
define the {\em restriction} of $J$ to $A$ to be the identity $I =
\langle A , F \cap ([A]^{2} \times [A]^{2}) \rangle$. This will be
written: $ I = J \rest A$.  Similarly for V-identities.  

\begin{defi}
Let $n < \omega$, $I=\langle A,E \rangle $, $J = \langle B,F\rangle$
be identities, $\bar{a} = \langle a_{1},\ldots, a _{n} \rangle \in
 \,^{n} A$ be a sequence of distinct elements from $A$, and
$\bar{b} = \langle b_{1},\ldots, b_{n} \rangle \in
 \,^{n} B$ be a sequence of distinct elements from $B$. $J$ is obtained
from $I$ by  {\em duplication} of $\bar{a}$ to $ \bar{b}$ if 
\begin{enumerate}
\item $B= A\sqcup \bar{b}$
\item $I = J\rest A$
\item the mapping which is the identity on $A\setminus \bar{a}$ and
  which maps $\bar{a}$ to $\bar{b}$ is is an embedding of $I$ into $J$ 
  as structures
\item $F$ is the least equivalence relation on $[B]^2$ consistent with
  i) -- iii).
\end{enumerate}
When $n =1$ we say that $J$ is a {\em one-point} duplication of $I$.  
\end{defi} 

\begin{defi}
Let $I=\langle A,E,<_{A} \rangle$ and $J= \langle B,F,<_{B} \rangle$
be V-identities.  $J$ is obtained from $I$ by end-duplication if there
exist final segments $\bar{a}$, $\bar{b}$ of $\langle A,<_A\rangle$,
$\langle B,<_B\rangle$ respectively such that the structure $\langle B,F
\rangle$ is obtained from the structure $\langle A,E\rangle$ by
duplicating $\bar{a}$ to $\bar{b}$. 
\end{defi}

Let IDE$_{V}$ denote the minimal collection of V-identities which is
closed under end-duplication, the taking of subidentities, and which
contains the trivial V-identity of size one. Let IDE $ =\{ \langle
A,E\rangle : \langle A,E,< \rangle \in IDE_{V} \}$

Let IDM denote the least class of identities which is closed
under duplication, the taking of subidentities and which contains the
identity of size one.  We quote the following results of Shelah:

\begin{theo}[\cite{Shappendix}] $ {\cal I}(\aleph _{1} ) \supset$ IDE. 
\end{theo}

\begin{theo}[\cite{ShmodelsII}]
\labl{Alan} There exists $f:[\aleph _{1}]^{2} \longrightarrow \aleph _{0}$ such that ${\cal I} (f) \subseteq$ IDE.
\end{theo}
 
\begin{coro}\labl{theo2.3}
${\cal I}(\aleph _{1} ) =$ IDE.
\end{coro}
\begin{theo}\labl{shel}[\cite{Shtwocard}]
If $\kappa \geq  \aleph _{\omega}$ then $ {\cal I}(\kappa ) \supseteq$ IDM.
\end{theo}

\section{ The Partial Order}
	
	We are going to define a partial order that will be used as a forcing notion to allow us to omit an identity from ${\cal I}(\aleph _{m})$.  The construction of the partial order uses historical forcing \cite {ShStHist}.  In this method conditions are allowed into the partial order if they can be constructed from the amalgamation of simpler conditions satisfying certain properties. As a preliminary we construct a model which has many definable subsets.

	In the following, first-order languages will be denoted by $L$, possibly with subscripts.  The variables of each first-order language are the same and have a canonical well ordering of order type $\omega $.  We denote by $ F _{n} (L) $ the set of all formulas of $L$ in which at most the first $n$ variables occur free.  We let $F(L)= \bigcup \{ F _{n} (L): n < \omega \}$.  Let $L_{0} $ denote the language $\{ < \} $.  Let $< _{m}$ denote the natural ordering of $\aleph _{m} $.  We regard $\langle \aleph _{m} ; < _{m} \rangle $ as an $ L_{0}$-structure.  From now on fix $m$ to be an integer $ \geq 1$.
\begin {lemm}
	There exists a countable relational language $L \supset L_{0} $ and an expansion $ {\cal M } $ of $\langle \aleph _{m} ; < _{m} \rangle $ to $L$ such that for each $ k < \omega$ and $\phi \in F _{k+1} (L) $ there exists $\psi \in F _{k+2} (L) $ such that for every $k$-tuple $\bar{a} \in \aleph _{m}, $ the relation defined by $\psi ( x,y,\bar{a}) $ in ${\cal M} $ well orders $ \{ b \in {\cal M } : {\cal M } \models \phi ( b,\bar{a}) \} $ with order type $| \{ b \in { \cal M } : {\cal M } \models \phi ( b,\bar{a} ) \} |$.
\end{lemm}
\begin{proo} We define a sequence of models and languages $\langle \langle {\cal M}_{i}, L _{i} \rangle : i< \omega \rangle$ with the property:

	For each $k< \omega, \varphi \in F_{k+1} (L_{n})$ and $k$-tuple $ \bar{a }\in \aleph _{m}^{k}$ there exists $\psi \in F_{k+2} (L_{n+1}) $ such that $\psi ( x,y, \bar{a})$ well orders $\{ b : {\cal M} \models \varphi ( b,\bar{a})\}$ in order type $|\{ b : {\cal M} \models \varphi ( b,\bar{a})\}|$ .

	Assume that $\langle {\cal M}_{n},L_{n} \rangle$ has been defined.  For every $k<\omega $ and $\varphi \in F_{k+1} (L_{n}) $ we add a $(k+2)$-ary relation symbol $W_{\varphi}$ to $L_{n+1}$ and expand the structure ${\cal M}_{n} $ to ${\cal M}_{n+1} $ so that for all $\bar{a} \in \aleph _{m}^{k},$ in ${\cal M}_{n+1} $ the formula $W_{\varphi} ( x,y,\bar{a}) \mbox { well orders }  \{ b : {\cal M}_{n} \models \varphi (b,\bar{a} ) \} \mbox { in order type } |\{ b : {\cal M}_{n} \models \varphi (b,\bar{a} ) \} |.$  We let $L = \bigcup \{L_{n}: n < \omega \} $ and ${\cal M } = \bigcup \{ {\cal M}_{n}: n < \omega \}$.
\end{proo}

\begin{defi}
Let $b \in \aleph _{m}, k,n <\omega$. and $\bar{a} \in  \aleph _{m}^{k}.$  We say $ b R_{n}  \bar{a}  $ if there exists $\varphi \in  F(L)  $ such that $  {\cal M } \models (\exists ^{\leq \aleph _{n} }x )\varphi ( x,\bar{a})$ and ${\cal M } \models \varphi ( b,\bar{a})$. 
\end{defi}
\medskip

Let ${\Bbb R} = \{ \langle u,c \rangle : u \in [ \aleph _{m} ]^ {< \aleph _{0} } , c : [u] ^{2} \longrightarrow \omega \} .$  

\begin{defi}
$p = \langle u,c \rangle  \in { \Bbb R } $ is the {\em  amalgam} of $p^{0}= \langle u^{0} ,c ^{0} \rangle$ and $p^{1}= \langle u^{1} ,c ^{1} \rangle  \in {\Bbb R}$ if there exist $ h < \omega $ and increasing sequences $ i^{0} _{0} ,\ldots, i^{0} _{h}$ and $ i^{1} _{0}, \ldots, i^{1} _{h} $ in $ \aleph _{m} $  such that for all $s, t$
with $( 0 \leq s < t \leq h ),$ and all $i,j,k,l < \aleph _{m} $:
\begin{enumerate}
\item$ u^{0} = \{i^{0} _{0} \ldots i^{0} _{h} \}$ and $  u^{1} = \{i^{1} _{0} \ldots i^{1} _{h} \}$
\item $ c^{0}(\{ i^{0} _{s}, i^{0} _ {t}\}) = c^{1}(\{ i^{1} _ {s}, i^{1} _{t}\}) $
\item $ i^{0} _{t} = i^{1} _{t} \vee i^{0} _{t} < i^{1}_{t}$
\item $i^{0} _{t} \neq i^{1} _{ t}$ implies $ \neg i^{1} _{t} R_{0} u^{0}$
\item $u= u^{0} \cup u^{1} $
\item $ c \supset (c^{0} \cup c^{1})$
\item $(\{i,j\} \notin [u^{0}]^{2} \cup [u^{1}] ^{2})$ implies $ c(\{i,j\} ) \not\in  \rng (c^{0}) \cup \rng (c^{1})$ 
\item$ c(\{ i,j \})= c(\{ k,l \})$ implies $ (\{ i,j \}= \{ k,l \}  \vee \{ i,j \}, \{ k,l \} \in [u^{0}]^{2} \cup [u^{1}] ^{2})$

\end{enumerate}
\end{defi}

\begin{defi}
$q = \langle u ^{q} ,c ^{q} \rangle \in { \Bbb R } $ is a {\em
 one-point extension} of $p = \langle u ^{p},c^{p} \rangle \in {\Bbb
 R}$ if $ u^{q} = u^{p} \cup \{ r\}$ for some $r > u^{p}, \, c^{p}
 \subset c ^{q}$, and for all $ i,j,k,l \in u^{q}$
\begin{enumerate}
\item $ \{ i,j \} \not\in$ dom$(c^{p}) $ implies $ c^{q}(\{ i,j \})
\not\in$ rng$(c^{p})$ 
\item $  c^{q}(\{ i,j \}) = c^{q}( \{k,l\})$ implies $ (\{ i,j \}, \{
k,l \} \in$ dom$(p) \vee \{ i,j \} = \{ k,l\})$.
\end{enumerate}
\end{defi}

\begin{defi}
We now define a sequence of subsets of ${\Bbb R}$.   Let $ {\Bbb P} _{0} = \{ \langle u,c \rangle \in {\Bbb R } : |u| =1\}.$  Given ${\Bbb P} _{n}$ we let ${\Bbb P} _{n+1}$ be the subset of $ {\Bbb R} $  which contains ${\Bbb P} _{n}$, all amalgam of pairs of elements from ${\Bbb P} _{n}$ and all one-point extensions of elements of $ {\Bbb P} _{n}$.  Let ${\Bbb P } = \bigcup \{{\Bbb P } _{n}: n < \omega \}$.  Given $ p= \langle  u^{p},c^{p} \rangle $ and $ q= \langle u^{q}, c^{q} \rangle$ we let $ p \leq q$ if and only if $ u^{p} \supseteq u^{q}$ and $c ^{p} \supseteq c^{q}$.
\end{defi} 
  
\medskip

	It should be noted that the order of $p^{0}$ and $p^{1}$ in the definition of amalgamation is important because of the asymmetry in the properties of $p^{0}$ and $p^{1} $ required by the definition. It is also worth observing that in terms of the notion of duplication, the amalgam of $p^{0}$ and $p^{1} \in {\Bbb R}$ may be regarded as being obtained from $ p^{0} $ by simultaneous duplication of all the elements in $ u^{0} \setminus u^{1} $.  Of course, the amalgamation also requires that none of the elements of $ u^{0} \setminus u^{1}$ belong to the countable set definable over $ u^{0}$. The closure of ${\Bbb P} $ under one-point extensions is necessary to show that our forcing produces a function whose domain is of size $\aleph _{m} $.

Let $ p = \langle u,c \rangle \in {\Bbb P } $ and $ I = \langle A,E \rangle$
be an identity.  We say that $p$ {\em realizes } $I$ if the $ \omega$-coloring
$ c : [ u ] ^{2} \longrightarrow \omega $ realizes $I$.  The mapping
$h: A \longrightarrow u $ demonstrating the realization is called an
{\em embedding} of $I$ inot $p$.

\begin{lemm} 
${\Bbb P}$ is c.c.c
\end {lemm}
\begin{proo} Let $ \langle p_{\alpha} : \alpha < \omega _{1} \rangle $ be a sequence of conditions.  By thinning we can suppose that there are $n,l< \omega $ and $i^{\alpha } _{j} \; (\alpha < \omega _{1}, 0\leq j\leq n) $ such that for all $ \alpha , \beta < \omega _{1} $ and all $j,k $ with $ 0\leq j < k \leq n$

\begin{enumerate}
\item $u^{p_{\alpha }}= \{ i^{\alpha } _{0}, \ldots ,i^{\alpha} _{n} \} $
\item $ i^{\alpha } _{j} < i^{\alpha} _{k} $
\item $c^{p_{\alpha}} ( \{ i^{\alpha} _{j} , i^{\alpha} _{k} \} ) = c^{p_{\beta} }( \{ i^{\beta }_{j} ,i^{\beta} _{k} \})$ 
\item $ p_{\alpha} \in {\Bbb P} _{l}.$ 
\end{enumerate}
Applying the  $\Delta $-system argument allows us to thin the sequence of conditions further so that $$ \forall t( 0 \leq t \leq n \Rightarrow
[\forall \alpha \forall \beta ( i^{\alpha } _{t} = i^{\beta } _{t}) \vee( \forall \beta < \omega _{1})( \forall \alpha < \beta )   ( i^{\alpha } _{t} < i^{\beta } _{t})]). $$

Let $T= \{ t\leq n : i^{\alpha } _{t} \neq i^{\beta } _{t} $ some $\alpha ,\beta < \omega _{1} \}$.
	To prove that ${\Bbb P} $ is c.c.c. it is sufficient to show that $p_{0}$ and $p_{\alpha } $ are compatible for some $ \alpha < \omega _{1} $.  For $p_{0}$ and $p_{ \alpha } $ to have a common extension by definition we need only that
$ t \in T $ implies $ \neg i^{\alpha } _{t} R _{0} ( i^{0} _{0} \ldots i^{0} _{n} )$. Since the language is countable $ |\{ i \in \aleph _{m} : i R _{0} (i^{0} _{0} \ldots i^{0} _{n} ) \} | = \aleph _{0} .$  For each $t \in T , i^{\alpha } _{t}$ is strictly increasing in $\alpha $.  Hence $ \neg i^{\alpha } _{t}   R _{0} (i^{0} _{0} \ldots i^{0} _{n} )$ for all sufficiently large $ \alpha < \omega _{1}. $  Since $T$ is finite, the condition above is satisfied for all sufficiently large $\alpha $ .  This proves the lemma.
\end{proo}
\begin{lemm}
For each $\alpha < \aleph_{m}$  $$\{ \langle u,c \rangle\in {\Bbb P} : \exists \beta ( \beta \in u \wedge \alpha < \beta) \}$$ is dense in  $({\Bbb P},< )$.
\end{lemm}
\begin{proo} Clear from the definition of one-point extension.
\end{proo}

\begin{lemm}
\labl{t1}
Let $M$ be any model and $G$ be ${\Bbb P} $ generic.  In $M[G]$ there is a function $ f : [\aleph _{m}] ^{2} \longrightarrow
\omega $ such that every identity realized by $f$ is a subidentity of an identity realized by some $ p \in {\Bbb P}$.
\end{lemm}
\begin{proo} Standard forcing technique. 
\end{proo}
\section {Omitting $I_{m}$ }

Throughout this section let $m$ be an integer $\geq 1$.  We denote by
$ I_{m}$ the identity $\langle A,E \rangle $ where $ A = \, ^ { m+1}2$ and
$\{ \eta, \nu \} \simeq_{E} \{ \alpha ,\beta\}$ if and only if 
$ \eta \cap \nu = \alpha \cap \beta$. We show that the model of ZFC
produced by forcing with the partial order of the previous section contains
a function $ f : [ \aleph _{m+1} ] ^{2} \longrightarrow \omega $ that does not
realize $ I_{m}$.  It should be noted that $ I_{m} \in \IDM$. Since theorem \ref{shel} says that $ { \cal I } (\aleph _{ \omega } ) \subseteq \IDM$, the model produced by forcing clearly satisfies $ {\cal I } ( \aleph _{m+1})
\subsetneqq { \cal I } ( \aleph _{ \omega })$.

\begin{defi}
Let $ I= \langle A,E \rangle$ be an identity and $ p=\langle u,c \rangle \in {\Bbb R}$.  Then $ h:A \longrightarrow u $ is an { \em embedding of $I$ into $p$} if $h$ is an embedding of $I$ into the $ \omega$-coloring $c$.
\end{defi}
\begin{defi}
$ \langle \eta _{0}, \ldots , \eta _{m} , \eta _{m+1} \rangle $ is a {\em special sequence} if 
\begin{enumerate}
\item $\eta _{i} \in  \, ^{m+1}2$
\item $ | \eta _{i} \cap \eta _{i+1}| = i$ for all $ i\leq m.$
\end {enumerate}
\end{defi}
\begin {lemm}
\labl{pp}
	Let $p\in {\Bbb R}$ and $ h$ be an embedding of $I_{m}$ into $p$. There exists a special sequence $ \langle \eta _{0}, \ldots , \eta _{m} , \eta _{m+1} \rangle $ so that $h(\eta _{i}) R_{0}(h(\eta _{0}) \ldots h(\eta _{m-1}))$ for $ i \in \{ m,m+1 \}.$ 
\end{lemm}
\begin{proo}{proo}  We define $ \eta _{k} \in 2 ^{m+1} , 0 \leq k \leq m-1,$ by induction on $k$ such that
\begin {enumerate}
\item $| \eta _{i} \cap \eta _{i+1 } | = i $ for all $ i < m-1$
\item $h( \gamma ) R _{m-(i+1)} ( h( \eta _{0}), \ldots , h( \eta _{i}))$ for all $ i< m$ and all $ \gamma \in  \, ^{m+1}2$ such that $ | \gamma \cap \eta _{i} |=i.$
\end{enumerate}

Let $ \eta _{0} $ be the unique $ \nu \in  \, ^{m+1} 2 $ such that $ h ( \nu ) = $ max(rng($h$)).  Suppose that $ \eta _{k} $ has been suitably defined for all $ k \leq j$ where $ j < m-1$.  Let $C$ denote  $\{ \nu \in  \,^{m+1}2 : | \nu \cap \eta _{j} | =j \}$.  From the induction hypothesis there exists $ D \subset \aleph _{m} $ such that rng$( h \rest C ) \subset D$ and $ D $ definable in ${\cal M}$ over $\{ h ( \eta _{0}), \ldots , h( \eta _{j}) \}.$  From the choice of ${\cal M} $ there is a relation $ < _{D}$ definable over $\{ h ( \eta _{0}), \ldots , h( \eta _{j}) \}$ which well orders $ D $ in the order type less than or equal to $ \aleph _{m-(j+1)}.$  Let $ \eta _{j+1} $ be the unique $ \nu \in C $ such that $ h(\nu ) $ is the $<_{D}$ maximal element of rng$(h \rest C ).$  Clearly, $ | \eta _{j+1} \cap \eta _{j} | =j.$  Consider $ \gamma \in  \,^{m+1}2 $ such that $ | \gamma \cap \eta _{j+1} | = j+1.$  Clearly, $ \gamma \in C $ and $ \gamma \neq \eta _{j+1}.$  Hence $ h(\gamma ) \in D $ and $h( \gamma ) < _{D } h( \eta _{j+1} ).$  It follows that $ h( \gamma ) R _{m-(j+2)} ( h( \eta _{0}) , \ldots , h(\eta _{j}),h( \eta _{j+1} )).$  This completes the induction step and the definition of $ \eta _{0} , \ldots , \eta _{m-1}.$  Letting $ \eta _{m} , \eta _{m+1} $ be the two elements of $ \{ \nu \in  \, ^{m+1}2 : | \nu \cap \eta _{m-1} | =m-1 \} $ completes the proof.
\end{proo}
\begin{lemm}
\labl{qq}
Let $ \langle \eta _{0} , \ldots , \eta _{m+1} \rangle $ be a special sequence, $ p,q \in {\Bbb R} , p $ be a one-point extension of $q$, and $h$ be an embedding of $ J = I_{m} \rest \{ \eta _{0} , \ldots , \eta _{m+1} \}$ in $p$.  Then $h$ is an embedding of $J$ in $q$.  
\end {lemm}
\begin{proo}  Let $ p = \langle u,c \rangle $ and $ q = \langle v,d \rangle.$  Towards a contradiction suppose that $ u\setminus v = \{ h( \eta_{i} )\}. $  If $i<m $, then $ \{ \eta _{m}, \eta _{i} \} , \{ \eta _{m+1}, \eta _{i} \} $ get the same color in $ I_{m} $, but $c( \{ h( \eta _{m}) , h( \eta _{i} ) \} ) \neq c( \{ h( \eta _{m+1}) , h( \eta _{i} ) \} )$ since $p$ is a one-point extension of $q$.  This contradicts $h$ being an embedding.  If $i\in \{m,m+1\}$, the consideration of the pairs $ \{ \eta _{0}, \eta _{m} \}, \{ \eta _{0}, \eta _{m+1} \} $ leads to a similar contradiction.
\end{proo}

\begin {lemm}
\labl {t2}Let $p \in {\Bbb P}$.  Then there does not exist an embedding $h$ of $ I_{m}$ into $p$.

\end {lemm}

\begin{proo}  Towards a contradiction suppose the theorem fails.  From lemma~\ref{pp} there exist a special sequence $ \langle \eta _{0} , \ldots , \eta _{m+1} \rangle , p= \langle u,c \rangle \in { \Bbb P }, $ and $h$ embedding $ J = I_{m} \rest \{ \eta _{0} , \ldots , \eta _{m+1 } \}$ in $p$ such that $ h( \eta _{i} ) R _{0}  (h( \eta _{0}), \ldots , h( \eta _{m-1}))$ for $ i \in \{ m,m+1\}.  $  Fixing $ \langle \eta _{0} , \ldots , \eta _{m+1 } \rangle $ and $h$ choose $p$ to minimize $ |u|.$  From lemma~\ref{qq}, $p$ is not a one-point extension of $ q \in { \Bbb P}.$  Therefore there are $ p^{0} = \langle u^{0} , c^{0} \rangle , p^{1} = \langle u^{1} ,c^{1} \rangle \in { \Bbb P } $ such that $p$ is the amalgam of $ p^{0}$ and $ p^{1} $.  Since neither $ p^{0}$ nor $ p^{1} $ can replace $p$, there exist $ i,j \leq m+1$ and $ a,b \in u $ such that $ h( \eta _{i} ) = a \in u^{0} \setminus u^{1}, h( \eta _{j} ) = b \in u^{1} \setminus u^{0}$. 

From the definition of amalgamation, $ \{ a,b \} $ is the only pair in $ [u]^ {2}$ which is assigned the color $ c( \{ a,b \} ) $ by $p$.  The only pair in $[ \{ \eta _{0} , \ldots , \eta _{m+1} \} ] ^{2} $ which is assigned a unique color by $ I_{m} $ is $ \{ \eta _{m} , \eta _{m+1} \}.$  Without loss of generality $ i=m $ and $ j=m+1.$  Also, it is clear that $ h( \eta _{0}) , \ldots , h( \eta _{m-1} )$ are all in $ u^{0} \cap u^{1}.$

We have that $ b \in u^{1} \setminus u^{0} $ belongs to a countable set definable in $ { \cal M } $ over $ u^{0} \cap u^{1} $.  This contradicts the definition of amalgamation and completes the proof of the lemma.
\end{proo}
\begin{theo}
\labl{s1}
Let $M$ be any model and $G$ be ${\Bbb P}$ generic.  Then  $M[G]$ satisfies, $I_{m} \not\in {\cal I} (\aleph _{m}).$
\end {theo}
\begin{proo}  The result follows from  lemma \ref{t1} and lemma \ref {t2}. 
\end{proo}

\section{Inducing $ I_{m}$ }

To show that the identities induced at different cardinals can distinguish the cardinals themselves, as promised in the introduction, we will show that $I_{m} \in { \cal I} ( \aleph _{ m+1})$ for all $m$.  This will be done by showing that if a given collection of identities is induced at $ \aleph _{m}$ we can extend the collection in a nontrivial way and be assured that this new collection is induced at $ \aleph _{m+1}$.

\begin{defi}
Let $\langle J_{i}: 1 \leq i \leq n \rangle$ be a finite sequence of identities.  We define the {\em end-homogeneous amalgam} of the sequence as follows.  Choose a sequence of $ \omega$-colorings
$ c_{i} : [G _{i}]^{2} \longrightarrow \omega $ such that $  c_{i} \in J_{i}$, $ G_{i} \cap G _{j} = \emptyset $ for $ 1\leq i < j \leq n$, and $ \rng (c_{i}) \cap \rng (c_{j})= \emptyset $ for all $ 1 \leq i < j \leq n$.  
Let $ G = \bigcup \{ G _{i} : 1 \leq i \leq n \}$. Now choose a new 
$\omega$-coloring $  c: [G]^{2} \longrightarrow \omega $ such that
for all $ \{ r,s \}, \, \{ t,v \} \in [G] ^{2} $ and all $i, \, 1\leq i \leq n$,
\begin{enumerate} 
\item $ c \supset c_{i}$ 
\item $ c(\{ r,s \})  \in  $ rng$ (c_{i}) \text{ if and only if } \{ r,s \} \in $ dom$ (c_{i} )$  
\item  if $ \{ r,s \}, \, \{ t,v \} $ are not in $ \bigcup \{ \dom ( c _{j})
: 1 \leq j \leq n \} ,$ then 
$$c( \{ r,s \})= c(\{ t,v \} )  \Leftrightarrow \text{ min} \{ j: r \in G_{j} \vee s \in G_{j} \} = \text{ min} \{ j: t \in G_{j} \vee v \in G_{j} \}.$$
\end{enumerate}
The end-homogeneous amalgam of $ \langle J_{i}: 1 \leq i \leq n \rangle$ is the identity realized by $ c$.
\end{defi}

Let $ {\cal I }$ be a collection of identities. Define the {\em closure} of $ { \cal I } $ (denoted cl$({\cal I })$) to be the collection of identites produced by forming all end-homogeneous amalgam of all finite sequences of identities in ${\cal I }$.

\begin{theo}
Let $ 0< m < \omega$.  If ${\cal I } \subseteq {\cal I } ( \aleph _{m})$ then { \em cl}$({ \cal I } ) \subseteq {\cal I }(\aleph _{m+1})$.
\end{theo}
\begin{proo}  Let $ f: [\aleph_{m+1}] ^ {2} \longrightarrow \omega$ and $ \langle J_{i} :1 \leq i\leq n \rangle $ be a sequence of identities in ${\cal I } $.  We will produce by recursion a sequence $\langle \langle A _{k}, B _{k}\rangle :0 \leq k \leq n \rangle $ such that:
\begin {enumerate}
\item $ f $ induces $ J_{i} $ on the set $ A _{i} $ for $ 1 \leq i \leq n$
\item $ B_{i} \supset B _{i+1}$ for $ 0 \leq i < n$
\item $ | B _{i} | = \aleph _{m+1 } $ for $0 \leq i \leq n$
\item $ A _{i+1} \subset B _{i} \setminus B _{i+1}$
\item $ f(\{ a_{1},b_{1} \} ) = f( \{ a_{2}, b_{2} \}) $ whenever there exist $ i,j \;( 1\leq i\leq j \leq n) $ such that $ \{ a_{1},a_{2} \} \subset A_{i}$ and $\{ b_{1} ,b_{2} \} \subset B _{j}$.
\end{enumerate}

Define $ B _{0}$ to be $ \aleph _{m+1} $ and $A _{0} $ to be empty.  By induction suppose that $ \langle A_{i}, B _{i} \rangle $ have been defined for $ i \leq k < n.$  Let $ C _{k}$ be the first $ \aleph_{m} $ elements of $ B _{k}$.  For each $b \in B _{k} \setminus C _{k}$ there exists a subset $D _{k}$ of $ C _{k} $ and 
 $ c _{b,k } < \omega $ such that $ | D _{k} | = \aleph _{m} $ and $ f (\{ b, x \})= c _{ b,k} $  for all $ x \in D _{k} $.  Now choose  a finite set $ A _{b} \subset D _{k}$ such that $f$ induces $I_{k+1} $ on $A_{b}$.  There are only $ \aleph _{m}$ finite subsets of $ C _{k}$ and a countable collection of possible values
for $ c_{b,k}$.  Thus we can choose $B _{k+1} \subset B _{k} \setminus C _{k}$ of cardinality $ \aleph _{m+1} $ and $ c _{k} < \omega $ such that $ A_{b_{1} } = A _{b _{2}} $ for all $ \{ b_{1} , b _{2} \} \subset B _{k+1}$ and $ c_{b,k} = c _{k} $ for all $ b \in B _{k+1}.$ We let $ A _{k+1} = A _{b } $ for $ b \in B _{k+1}.$  It is easy to see that $f$ induces the desired identity on the set $ \bigcup \{ A _{i}: 1 \leq i \leq n \}$.  
\end{proo}

\begin{theo}
\labl{s2}
For all $m$ such that $ 1 \leq m < \omega$ we have $I_{m} \in {\cal I } ( \aleph _{m+1} ) $.
\end{theo}
\begin{proo}  The proof is by induction on $ k< \omega $. To show that $ I_{1} \in {\cal I } ( \aleph _{2} ) $ we let $ J_{1} =J_{2} $ be the trivial identity on two points and produce the identity J, the end-homogeneous amalgam of the sequence $ \langle J_{1} ,J_{2} \rangle$.  An analysis of $J$ shows it to be $I_{1}$.  By induction we assume that $I_{k} \in {\cal I} ( \aleph _{k+1}).$  We then let $ J_{1} = J_{2} = I_{k} $ and produce an identity $K$, the end-homogeneous amalgam of the sequence $ \langle J_{1} , J_{2} \rangle.$  By the theorem we have that $ K \in {\cal I } ( \aleph _{k+2} ).$  An analysis of $K$ shows it to be $ I_{k+1}$.  
\end{proo}

The following theorem follow from theorem \ref{s1} and theorem \ref{s2}.

\begin{theo}
The consistency of $ZFC$ implies the consistency of $ {\cal I }( \aleph _{m}) 
\subsetneq {\cal I}( \aleph _{m+1})$ for $ m \geq 1$.
\end{theo}

\medskip
The question that now remains to be answered is whether or not it is 
consistent that $ {\cal I } (\aleph _{2}) = {\cal I}( \aleph _{\omega})$.  
Since this is trivially true when $CH$ holds, the question is only valid when 
a large continuum is also demanded. 
We have been successful in getting the consistency of the above statement when the continuum is large and the result will appear in the future.

\end{document}